\documentclass[a4paper,10pt]{amsart}
\usepackage[english]{babel}

\usepackage[curve]{xy}
\splinetolerance{1pt}

\newtheorem{definition}{Definition}
\newtheorem{proposition}{Proposition}

\newtheorem{lemma}{Lemma}
\newtheorem{corollary}{Corollary}
\newtheorem{theorem}{Theorem}

\newcommand{\un}{\mathbf{1}}
\newcommand{\zero}{\mathbf{0}}
\newcommand{\To}{\rightarrow}
\newcommand{\supp}[1]{{\it supp}(#1)}

\newcommand{\Co}{\mathbb{C}}
\newcommand{\Proj}{\mathbb{P}}
\newcommand{\Cot}{\mathbb{C}^{\times}}
\newcommand{\Z}{\mathbb{Z}}
\newcommand{\Cox}{\Co[x,\frac 1 x]}
\newcommand{\CoxS}{\Co_{s}[x,\frac 1 x]}
\newcommand{\eqmod}[3]{#1 \equiv #2 (\emph{mod}\;#3 )}

\newcommand{\centr}[1]{Z\left(#1\right)}

\newcommand{\rend}[2]{\mathcal{E}nd\left( #1 \right)_{#2}}

\newcommand{\matr}[2]{\mathcal M_{#1}\left(#2\right)}
\newcommand{\kernel}[1]{\mathcal{K}er{#1}}
\newcommand{\image}[1]{\mathcal{I}m{#1}}
\newcommand{\fact}[2]{{#1}\Big/{#2}}

\newcommand{\TL}{\bf{T}_N(\delta)}
\newcommand{\TLN}[1]{{\bf T}_{#1}(\delta)}
\newcommand{\varietyN}{\mathcal{V}_N}
\newcommand{\variety}[1]{\mathcal{V}_{#1}}
\newcommand{\varietyI}{\mathcal{W}_N}
\newcommand{\subalgN}{\mathcal{A}_N}
\newcommand{\idealN}{\mathcal{I}_N}
\newcommand{\overN}{\mathcal{K}_N}
\newcommand{\posetN}{\mathcal T_N}
\newcommand{\poset}[1]{\mathcal T_{#1}}
\newcommand{\diag}[2]{D(#1 \rightarrow #2)}
\newcommand{\cut}{{\it l}}

\title[Description of the center of the affine Temperley-Lieb
algebra]{Description of the center of the affine Temperley-Lieb algebra of
type $\widetilde{A_N}$}
\author{Masha Vlasenko}
\email{mariyka@imath.kiev.ua}
\address{Institute of Mathematics of NAS of Ukraine, Tereshchenkivska str. 3,
Kyiv, 01601, Ukraine}

\begin{document}

\begin{abstract}
Construction of the diagrammatic version of the affine Temperley-Lieb
algebra of type $\widetilde{A_N}$ as a subring of matrices over the Laurent polynomials is given. We move towards geometrical
understanding of cellular structure of the Temperley-Lieb algebra.
We represent its center as a coordinate ring of the certain affine 
algebraic variety and describe this variety constructing its
desingularization.
\end{abstract}

\maketitle

\section*{Introduction}
In this paper we deal with a sequence of infinite-dimensional algebras $\TL$, $N
= 1,2,3,\dots$, introduced by Graham and Lehrer in \cite{GrahamLehrer1998}.
These algebras are extended versions of the Temperley-Lieb quotients of the
affine Hecke algebras of type $\widetilde{A_N}$. Let us recall the
naive definition of $\TL$.

Consider a vertical cylinder with $n$ marked points on the top circle
component of its boundary and $m$ marked points on the bottom circle
component, with $\eqmod nm2$. If we join these marked points pairwise by arcs on
the (lateral) surface of the cylinder without intersection we get an
\emph{affine diagram}. One can also add to affine diagram a finite number of
circles which circumnavigate the cylinder, if it is possible to do without
intersections with arcs. Let us denote by $\diag mn$ the set of homotopy
classes of such an affine diagrams. Take $\alpha \in \diag mn$, $\beta \in \diag
km$. If we glue the bottom boundary circle of $\alpha$ to the top boundary
circle of $\beta$ so that corresponding marked points coincide, we get a
cylinder again, and there can be some loops on its surface. Let us denote the
number of loops by $m(\alpha,\beta)$, and the affine diagram obtained after
removing all the loops by $\alpha \circ \beta$. Then $\alpha \circ \beta \in
\diag kn$.

\begin{definition}
Let $\delta \in \Co$ be a complex number. $\TL$ is an associative $\Co$-algebra
spanned over $\Co$ by affine diagrams $\diag NN$ with multiplication
$$
\alpha \beta = \delta^{m(\alpha,\beta)} \alpha \circ \beta
$$
for $\alpha,\beta \in \diag NN$.
\end{definition}

One can find the proof that this multiplication really agrees with associativity
and diagrams $\diag NN$ are linearly independent in \cite{GrahamLehrer1998}.
Note that there is a unity $\un \in \TL$ given by diagram $\un_N$ on the picture
in section~1.

Let us enumerate the marked points on the top and bottom circles of cylinder by
$1,\dots,N$ sequently, so that top $1$ is over bottom $1$, and so on. Take
some point between $N$ and $1$ on the top and some point between $N$ and $1$
on the bottom circle, consider a line segment $l$ on the surface
connecting these points. For the diagram $\alpha \in \diag NN$ its \emph{rank}
$|\alpha|$ is a minimal number of intersections with $l$ among all the
diagrams from the homotopy class $\alpha$. A diagram is called \emph{monic} if
it has no arcs connecting bottom point to bottom point. It was proved in
\cite{FanGreen1999} that the Temperley-Lieb quotient of the affine Hecke algebra
of type $\widetilde{A_N}$ is isomorphic to a subalgebra of $\TL$, spanned by all
non-monic diagrams of even rank together with the unity $\un_N$.

In section~1 of this paper we realize $\TL$ as a subalgebra of matrices over the
ring of Laurent polynomials in one variable $\Cox$ and construct a large simply
arranged ideal inside of $\TL$ (Theorem~\ref{ideal_inside}).

Denote $\posetN = \{ t \in \Z \;|\; 0 \le t \le N, \eqmod{t}{N}2 \}$. In
section~2 we obtain the following result.

\begin{corollary}\label{variety_decomposition}
The center of $\TL$ is a coordinate ring of a certain affine algebraic variety.
This variety has $\left[\frac n2\right]+1$ irreducible components. Each irreducible component is a curve with finite number of singular points.

\end{corollary}

We also construct a certain desingularization of this variety
(Theorem~\ref{desingularization}), and propose an algorithm to calculate the
numbers of intersections (self-intersections) of its irreducible components. In
section~3 we perform these calculations for $\TLN 3$.

Representation theory of algebras $\TL$ is known: in
\cite{GrahamLehrer1998} the set of all irreducible modules is described,
and the composition factors of so called \emph{cell modules} are determined.
They use concept of \emph{cellular algebras} (\cite{GrahamLehrer1996}). If some
finite dimensional algebra has cellular structure its representation theory can
be obtained in standard way. The generalization of the concept for infinite
dimensional algebras is possible but requires some special technique in each
concrete case (\cite{GrahamLehrer200?}). We also use cellular structure of $\TL$
in our constructions (section~1), but we move towards its geometrical
description. The problem of geometrical understanding of cellular structure in
general was stated in \cite{Lehrer}.


\section{Realization of $\TL$ as a subalgebra of matrices over the ring
of Laurent polynomials in one variable $\Cox$}

For $k \in \posetN$ let $A_k$ be a linear span over $\Co$ of diagrams $\diag
kN$. Then $A_k$ is a left $\TL$-module with the action of diagrams on diagrams
by multiplication:
\begin{equation*}
\alpha \beta = \delta^{m(\alpha,\beta)} \alpha \circ \beta
\end{equation*}
for $\alpha \in \diag NN, \beta \in \diag kN$. On the other hand, $A_k$ is a
right $\TLN k$-module. Multiplication law is associative
(\cite{GrahamLehrer1998}), so these actions commute.

We call an arc a \emph{through arc} if it connects bottom point to top
point. Let $B_k \subset A_k$ be a linear span of diagrams having strictly
less then $k$ through arcs. Note that the number of through arcs in $\alpha
\circ \beta$ is not greater then their number in each of the diagrams $\alpha,
\beta$. Thus $B_k$ is preserved under both left $\TL$- and right $\TLN
k$-action. So $C_k = \fact {A_k}{B_k}$ is a
$\TL$~-~$\TLN k$-bimodule as $A_k$ and $B_k$ are ones.

Note that $C_k$ is spanned over $\Co$ by all monic diagrams from $\diag kN$. Now
we consider left action of $\TLN k$ on $C_k$.

Let $k>0$. Denote by $I^{k-2}_k$ a linear span of all non-monic diagrams
in $\TLN k$. $I^{k-2}_k$ is an ideal, and all its elements have zero action on
$C_k$. Consider the diagram $\tau_k$, usually
it is called \emph{twist}:

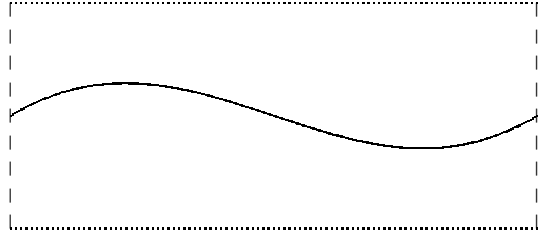
\[\begin{xy}
<0cm,1cm>;
<7cm,1cm>**@{.};
<7cm,4cm>**@{--};
<0cm,4cm>**@{.};
<0cm,1cm>**@{--},
(0,25);
(5,40)*{\bullet}**@{-},
(5,10)*{\bullet};
(15,40)*{\bullet}**@{-},
(15,10)*{\bullet};
(25,40)*{\bullet}**@{-},
(55,10)*{\bullet};
(65,40)*{\bullet}**@{-},
(65,10)*{\bullet};
(70,25)**@{-},
(5,6)*{1},
(15,6)*{2},
(55,6)*{k-1},
(65,6)*{k},
(5,43)*{1},
(15,43)*{2},
(25,43)*{3},
(65,43)*{k},
(33,0)*{\tau_k: k \rightarrow k \;(k>0)},
\end{xy}\]
\[
\begin{xy}
(0,10);
(70,10)**@{.};
(70,40)**@{--};
(0,40)**@{.};
(0,10)**@{--},
(0,25);
(70,25)**\crv{(25,40)&(45,10)},
(43,0)*{\tau_0: 0 \rightarrow 0},
\end{xy}\]
The quotient $\fact{\TLN K} {I^{k-2}_k}$ is isomorphic to the group
algebra of $\Z$, or the algebra of Laurent polynomials $\Cox$. Indeed, this
quotient is spanned by all monic diagrams $\diag kk$. The monic diagram $\tau_k$
is invertible: $\tau_k^{-1}$ is just the diagram obtained by reflection of
$\tau_k$ in a horizontal line. And there is no monic diagrams in $\diag kk$
except powers of $\tau_k$, powers of $\tau_k^{-1}$ and the unity $\un_{k}$:
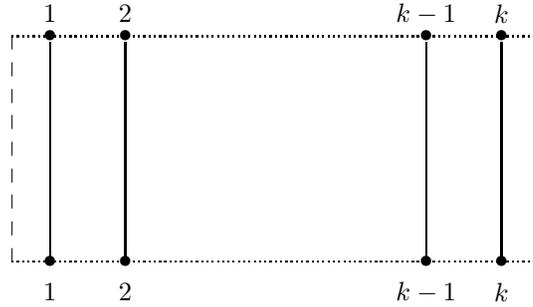
\[
\begin{xy}
<0cm,1cm>;
<7cm,1cm>**@{.};
<7cm,4cm>**@{--};
<0cm,4cm>**@{.};
<0cm,1cm>**@{--},
(5,10)*{\bullet};
(5,40)*{\bullet}**@{-},
(15,10)*{\bullet};
(15,40)*{\bullet}**@{-},
(55,10)*{\bullet};
(55,40)*{\bullet}**@{-},
(65,10)*{\bullet};
(65,40)*{\bullet}**@{-},
(5,6)*{1},
(15,6)*{2},
(55,6)*{k-1},
(65,6)*{k},
(5,43)*{1},
(15,43)*{2},
(55,43)*{k-1},
(65,43)*{k},
(33,0)*{\un_k: k \rightarrow k}
\end{xy}\]
Thus $C_{k}$ becomes a right $\Cox$-module with $x$ acting by multiplication by
$\tau_k$ from the right. It is a free right $\Cox$-module of finite
rank. Indeed, let us call a monic diagram \emph{standard} if it has no
circumnavigate circles and can be drawn without intersections of through arcs
with the line segment $\cut$. Now $k > 0$ and every monic diagram $\alpha \in
\diag kN$ can be uniquely "twisted" from the right to the standard one: there
exist a unique $r \in \Z$ such that $\alpha \tau_k^r$ is standard. Thus $C_{k}$
has a basis of standard diagrams over $\Cox$.

Let $k=0$. $\TLN 0$ itself is isomorphic to the algebra of polynomials in one
variable $\Co[x]$ if we define a clean diagram (without arcs and circles) to be
a unity $\un_0$. Indeed, $\diag 00$ consists of powers of $\tau_0$. So, $C_0$ is
a free right module over $\Co[x]$ with basis of standard diagrams.

To avoid cumbersome notations it would be useful to consider not $C_0$ but a
similar module over $\Cox$. Let $\widetilde{C_0}$ be a free right module over
$\Cox$ with basis of all standard diagrams from $\diag 0N$. We can define left
action of $\TL$ on $\widetilde{C_0}$. We prefer following symmetric version of
this action to be consistent with computations made in \cite{GrahamLehrer1998}:
let $\alpha \in \diag NN$, let $\beta \in \diag 0N$ is standard, and
$\alpha \circ \beta = \gamma \tau_0^r$ where $\gamma$ is also standard. Then put
$$
\alpha \beta = \delta^{m(\alpha,\beta)}\gamma (x + \frac 1x)^r.
$$
One can extend this to the action of $\TL$ on $\widetilde{C_0}$ by linearity.
Indeed, if $\varrho,\alpha \in \TL$ and $\beta \in \diag 0N$ is standard then
in $C_0$ one has
$$
\varrho (\alpha \beta) = (\varrho \alpha)\beta = \sum_{s}\gamma_s Q_s(x)
$$
for some standard $\gamma_s$ and $Q_s \in \Co[x]$. Then in $\widetilde{C_0}$
$$
\varrho (\alpha \beta) = (\varrho \alpha)\beta = \sum_{s}\gamma_s Q_s(x+\frac1x).
$$
So, $\widetilde{C_0}$ is a $\TL$~-~$\Cox$-bimodule, free over $\Cox$.

We denote
by $\CoxS \subset \Cox$ the subset of polynomials symmetric in $x$ and
$\frac1x$. Then $\CoxS \cong \Co[x]$.

Let us denote
$$
M_k = \begin{cases}
    C_k &\hbox{if $k>0$},\\
    \widetilde{C_0} &\hbox{if $k=0$}.\\
      \end{cases}
$$
Consider
$$
M = \prod_{k \in \posetN}{M_k}
$$
where the direct product is considered in the category of
$\TL$~-~$\Cox$~-~bimodules. Now we consider the action of $\TL$
on $M$, we claim it is faithful.

\begin{proposition}\label{faithful}
The kernel of the map $\TL \To \rend {M}{\Cox}$ is zero.
\end{proposition}
\begin{proof}
Every element in $b \in \TL$ is a unique linear combination of diagrams. Let us
denote the set of this diagrams by $\supp b$. Also define $\deg b$ to be a
maximal number of through arcs among $\supp b$, and a \emph{principal part}
$\hat{b}$ to be a linear combination of diagrams with $\deg b$ through arcs
such that $\deg (b-\hat{b}) < \deg b$.

Suppose there is a nonzero element $a \in \TL$ acting on $M$ as zero. Then
$\hat{a} \in \TL$ is nonzero. Let $k = \deg a$. Since $a - \hat{a}$ has
zero action on $M_k$, action of $\hat a$ on $M_k$ coincides with action of $a$
on $M_k$ and is zero. We are going to deduce that $\hat{a}$ is zero then. So
it will lead to a contradiction.

Let $\{ \beta^k_i; i = 1,2,\dots,d_k \}$ be the set of standard diagrams from
$\diag kN$. Then
$$
\hat{a} = \sum_{i,j=1}^{d_k}\beta^k_i P_{i,j}(\tau_k) (\beta^k_j)^*
$$
for some Laurent polynomials $P_{i,j}$, where for the diagram $\alpha \in
\diag mn$ we define $a^* \in \diag nm$ to be the reflection of $\alpha$ in a
horizontal line. Moreover, the polynomials $P_{i,j}$ are defined uniquely. Note
that $(\beta^k_j)^*\beta^k_s$ is a scalar multiple of a monic
diagram from $\diag kk$. We have identified the linear span of monic diagrams
from
$\diag kk$ with $\Cox$ (with $\CoxS$ for
$k=0$). Define the following $d_k \times d_k$ matrix over $\Cox$
\begin{equation}\label{mod_k_matrix}
R_k = \left( (\beta^k_j)^*\beta^k_s  \right)_{j,s=1}^{d_k} \; \in
\matr {d_k}{\Cox}
\end{equation}
in sense of this identification. Then the matrix of action of
$\hat a$ in basis $\beta^k_s$ is
\begin{equation}\label{action_matrix}
\begin{cases}
\left( P_{i,j}(x) \right)_{i,j=1}^{d_k} R_k &\hbox{for $k>0$},\\
\left( P_{i,j}(x+\frac 1x) \right)_{i,j=1}^{d_0} R_0 &\hbox{for $k=0$}.\\
\end{cases}
\end{equation}
This matrix should be zero. The calculation in
\cite{GrahamLehrer1998} (Corollary 4.7) shows that $\det R_k$
is a nonzero Laurent polynomial. Thus $P_{i,j} = 0$ for all
pairs $i,j$, and $\hat {a}$ is zero. This is a contradiction.
\end{proof}

So we get an inclusion
$$
\TL \subset \prod_{k \in \posetN}\matr{d_k}{\Cox},
$$
as $\rend M {\Cox} \cong \prod_{k \in \posetN}\matr{d_k}{\Cox}$,
where $d_k$ is a number of standard diagrams in $\diag kN$. This
number is obviously independent of choice of the line segment
$\cut$ of cut on the surface of the cylinder. It was calculated
in~\cite{GrahamLehrer1998}~(Proposition 1.11):
$$
d_k = \begin{pmatrix}N\\ \frac{N-k}2\\ \end{pmatrix}.
$$

We proceed to know how large is the image of this inclusion.

\begin{proposition}\label{centralizer}
The commutant of $\TL$ in $\prod_{k\in \posetN}\matr{d_k}{\Cox}$ is just
$\prod_{k \in \posetN} \Cox$, i.e. the set of matrices, scalar in each
matrix component.
\end{proposition}
\begin{proof}
Take an arbitrary matrix $A = \prod_{k \in \posetN}A_k$ commuting with
$\TL$. For an arbitrary $1 \le i,j \le d_k$ consider the element $b^{i,j} =
\beta^k_i(\beta^k_j)^* \in \TL$, let $B^{i,j} \in \matr {d_k}{\Cox}$ be the
matrix of action of $b^{i,j}$ on $M_k$. Then $A_k B^{i,j} = B^{i,j} A_k$. But $\image {B^{i,j}} = \beta^k_i \Cox$, so $A_k \beta^k_i = Q_i(x) \beta^k_i$ for
some $Q_i \in \Cox$. Recall nonzero monomials $(R_k)_{i,i} \in \Cox$ defined
by~(\ref{mod_k_matrix}). Then
$$
Q_i(x) (R_k)_{j,j}(x) \beta^k_i = Q_i(x) \beta^k_i (\beta^k_j)^*\beta^k_j = A_k
B^{i,j} \beta^k_j = B^{i,j} A_k \beta^k_j
$$
$$
= Q_j(x) B^{i,j} \beta^k_j = Q_j(x) (R_k)_{j,j}(x) \beta^k_i,
$$
hence $Q_i = Q_j$. So, $A_k$ is a scalar matrix.
\end{proof}

We often identify elements of $\TL$ with corresponding matrices over $\Cox$ in basis of standard diagrams.

Now we are going to construct a large ideal in $\TL$. For these purposes we
will define a set of Laurent polynomials $\{P_k(x); k \in \posetN \}
\subset \Cox$ such that for $k>0$ if $A \in  \matr {d_k}{\Cox}$ is divisible by
$P_k(x^k)$ then the matrix equal to $A$ in $k$-th component and zero in all other
components belongs to $\TL$. The $0$-th matrix component for each element of
$\TL$ is a matrix with symmetric polynomials in its entries, so $P_0$ will be a
polynomial of $x + \frac 1x$ such that any $d_0 \times d_0$ matrix $A$ of
polynomials of $x+\frac 1x$ divisible by $P_0$ will belong to $\TL$ when we
put zeros in all other matrix components.

We denote $\Cot = \Co \backslash \{0\}$. Consider $q \in \Cot$ such that $\delta
= -q-\frac 1q$. Of coarse every symmetric Laurent polynomial in $q$ is a
polynomial in $\delta$. Put
$$
G_k(x) = \prod_{r \in \posetN, r > k} (x^2 + \frac 1 {x^2}- q^r - \frac
1{q^r})^{d_r}
$$
for $k < N$ and $G_N(x)=1$.
For each $k \in \posetN$ we defined a matrix $R_k \in \matr{d_k}{\Cox}$ by
(\ref{mod_k_matrix}). It was calculated in \cite{GrahamLehrer1998}(Corollary
4.7) that either
\begin{equation}\label{det}
\det R_k = G_k(x) \hbox{ or } \det R_k = -G_k(x).
\end{equation}

For $k>0$ take $H_k \in \Cox$ to be the Laurent polynomial of smallest possible
degree such that $H_k(x^k)$ is divisible by $G_k(x)$, put $H_0(x)=G_0(x)$. Put
$$
P_k(x) = \prod_{r \in \posetN, r \le k} H_r(x).
$$

\begin{theorem}\label{ideal_inside}
There are inclusions: for odd $N$
$$
\prod_{k \in \posetN} P_k(x^k) \matr{d_k}{ \Cox } \subset \TL \subset
\prod_{k \in \posetN}\matr{d_k}{ \Cox },
$$
and for even $N$
$$
{P_0(x) \matr {d_0}{\CoxS}} \times {\prod_{k \in \posetN, k>0} P_k(x^k)
\matr{d_k}{ \Cox }} \subset \TL
$$
$$ \subset \matr {d_0}{\CoxS}\times
\prod_{k \in \posetN, k>0}\matr{d_k}{ \Cox }.
$$
\end{theorem}
\begin{proof}
Consider a filtration of $\TL$ by the set of ideals
$$
\TL = I^{N}_N \supset I^{N-2}_N \supset I^{N-4}_N \supset \dots
$$
where for $k \in \posetN$
$$
I^k_N = \left\{ a \in \TL | \deg a \le k \right\}.
$$

Consider a homomorphism $\pi_k: \TL \To \matr{d_k}{\Cox}$ where $\pi_k(a)$ is
the matrix of action of $a \in \TL$ on $M_k$ written in basis of standard
diagrams.

Note that if $k < s$ then $\pi_s(I^k_N) = \zero$. Recall the matrix $R_k \in
\matr{d_k}{\Cox}$ defined by (\ref{mod_k_matrix}). Then
$$
\pi_k(I^k_N) = \begin{cases}
        \matr {d_k}{\Cox} R_k &\hbox{for $k>0$},\\
        \matr {d_0}{\CoxS} R_0 &\hbox{for $k=0$}\\
               \end{cases}
$$
due to~(\ref{action_matrix}).

The theorem statement is a consequence of the following more general one:
if $A \in \matr{d_k}{\Cox}$ ($A \in \matr{d_0}{\CoxS}$) is devisible by
$P_k(x^k)$ (by $P_0(x)$) then there exist an element $a \in I^k_N$ with
$\pi_k(a) = A$ and $\pi_m(a)=\zero$ for $m \in \posetN$, $m \ne k$.
We are going to prove it by induction moving up in $k \in \posetN$.

The base of induction (either $k=0$ or $k=1$) follows from the description of
$\pi_s(I^{k})$ for $s \ge k$ given above. Indeed, for these cases $P_k(x)$
is devisible by $\det R_k$ due to~(\ref{det}). Consider the
matrix
\begin{equation}\label{inverse}
N_k = \left(N_{i,j} = (-1)^{i+j}\widehat{R_k^{j,i}} \right)_{i,j=1}^{d_k},
\end{equation}
where $\widehat{R_k^{i,j}}$ is a minor of matrix $R_k$ obtained by deleting of
$i$-th row and $j$-th column. Then $N_0 \in \matr{d_0}{\CoxS}$ and $N_k
\in \matr{d_k}{\Cox}$ for $k>0$, and $N_k R_k = \det R_k(x)$. So,
$$
\pi_1(I^1_N) = \matr{d_1}{\Cox}R_1 \supset \matr{d_1}{\Cox}N_1 R_1 \supset P_1(x)
\matr{d_1}{\Cox},
$$
$$
\pi_0(I^0_N) = \matr{d_0}{\CoxS}R_0 \supset \matr{d_0}{\CoxS}N_0 R_0 \supset
P_0(x)
\matr{d_0}{\CoxS}.
$$

Suppose we have proved the statement for all $k\in \posetN$, $k < r$. Take an
arbitrary matrix $B \in \matr {d_r} {\Cox}$ and consider $A = P_r(x^r)B$.
Take $Q(x) = \frac{H_r(x^r)}{\det R_r}$. It is a Laurent polynomial due to~(\ref{det}). Then $A = Q(x)P_{r-2}(x^r)B N_r R_r$, so $A
\in \pi_r(I^r_N)$. Let $C = Q(x)B N_r$. Recall that $\{ \beta^r_i;
i=1,2,\dots,d_r\}$ is a basis of standard diagrams in $M_r$. Consider element
$$
a_r = \sum_{i,j=1}^{d_r}\beta^r_i P_{r-2}(\tau^r_r) C_{i,j}(\tau_r)
(\beta^r_j)^*.
$$
Then $a_r \in I^r_N$ and $\pi_r(a_r) = P_{r-2}(x^r)C R_r = A$ due
to~(\ref{action_matrix}). So,
for $s > r$ we have $\pi_s(a_r)=0$ as $a_r \in I^r_N$. We claim that
$\pi_k(a_r)$ is devisible by $P_k(x^k)$ for $r>k>0$ and $\pi_0(a_r) = 0$ if
$N$ is even. This will finish the proof. Indeed, by inductional
assumption we can take $a_k \in I^{k}_N$ for each $k \in \posetN$,$r>k>0$ such
that $\pi_s(a_k)=\zero$ for all $s \ne k$ and $\pi_k(a_k)=-\pi_k(a_r)$. Then the
element $\sum_{0<k\le r}a_k$ will satisfy our conditions.

To prove the claim note that for arbitrary monic diagrams $\mu,\nu \in \diag
rN$ we have for $k \in \poset r$
$$
\pi_k \left( \mu \tau^r_r \nu^* \right) = x^k
\pi_k \left( \mu \nu^* \right).
$$
This statement is equivalent to $\mu \tau^r_r \nu^* \beta^k_t = \mu \nu^*
\beta^k_t \tau^k_k$ for $1 \le t \le d_k$. For $k=0$ one has $\mu \tau^r_r \nu^*
\beta^k_t = \mu \nu^* \beta^k_t$ obviously. For $k>0$ the substitution
of $\tau^r_r$ by $\tau_k^k$ is shown on the following picture:

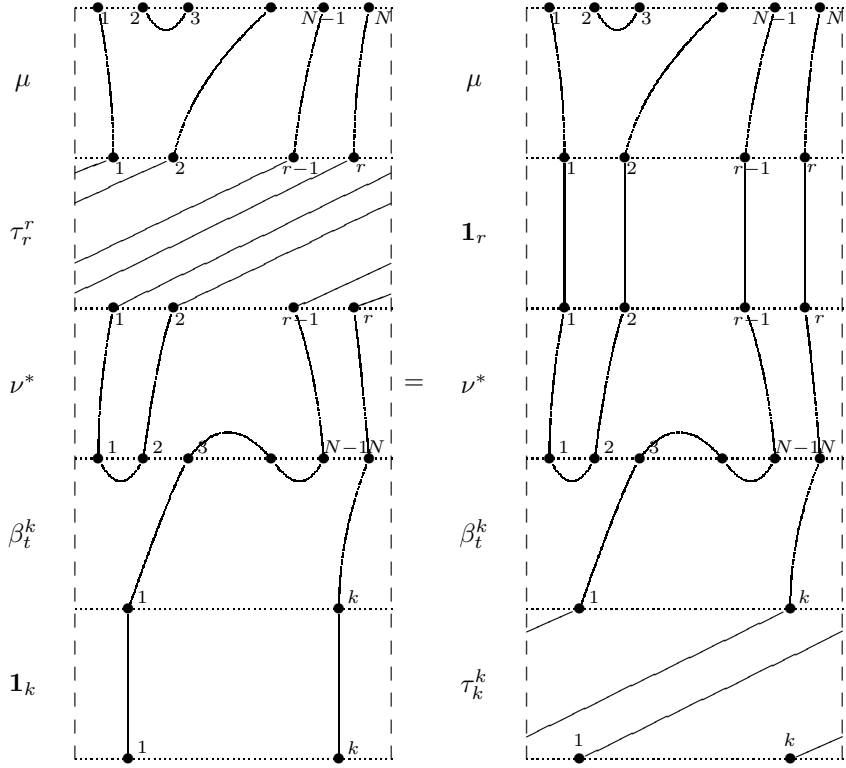
\[
\begin{xy}
(7,0);
(49,0)**@{.};
(7,20);
(49,20)**@{.};
(7,40);
(49,40)**@{.};
(7,60);
(49,60)**@{.};
(7,80);
(49,80)**@{.};
(7,100);
(49,100)**@{.};
(7,0);
(7,100)**@{--};
(49,0);
(49,100)**@{--};
(0,10)*{{\bf 1}_k},
(16,1)*{^1},
(14,0)*{\bullet};
(14,20)*{\bullet}**@{-},
(44,1)*{^k},
(42,0)*{\bullet};
(42,20)*{\bullet}**@{-},
(0,30)*{\beta^k_t},
(10,40)*{\bullet};
(16,40)*{\bullet}**\crv{(13,34)};
(16,21)*{^1},
(14,20);
(22,40)*{\bullet}**\crv{(15,35),(20,37)};
(33,40)*{\bullet};
(40,40)*{\bullet}**\crv{(37,34)};
(44,21)*{^k};
(42,20);
(46,40)*{\bullet}**\crv{(42,30)};
(0,50)*{\nu^*},
(12,41)*{^1},
(10,40);
(12,60)*{\bullet}**\crv{(10,50)};
(18,41)*{^2},
(16,40);
(20,60)*{\bullet}**\crv{(18,55)};
(24,41)*{^3},
(22,40);
(33,40)**\crv{(27,47)};
(43,41)*{^{N-1}};
(40,40);
(36,60)*{\bullet}**\crv{(39,52)};
(47,41)*{^N};
(46,40);
(44,60)*{\bullet}**\crv{(45,52)};
(0,70)*{\tau^r_r},
(7,78);
(12,80)*{\bullet}**@{-};
(13,78)*{^1};
(7,74);
(20,80)*{\bullet}**@{-};
(21,78)*{^2};
(7,66);
(36,80)*{\bullet}**@{-};
(37,78)*{^{r-1}};
(7,62);
(44,80)*{\bullet}**@{-};
(45,78)*{^r};
(49,78);
(12,60)**@{-};
(13,58)*{^1};
(49,74);
(20,60)**@{-};
(21,58)*{^2};
(49,66);
(36,60)**@{-};
(37,58)*{^{r-1}};
(49,62);
(44,60)**@{-};
(46,58)*{^r};
(0,90)*{\mu},
(12,80);
(10,100)*{\bullet}**\crv{(12,90)};
(11,98)*{^1};
(15,98)*{^2};
(16,100)*{\bullet};
(22,100)*{\bullet}**\crv{(19,94)};
(23,98)*{^3};
(20,80);
(33,100)*{\bullet}**\crv{(22,90)};
(36,80);
(40,100)*{\bullet}**\crv{(37,90)};
(40,98)*{^{N-1}};
(44,80);
(46,100)*{\bullet}**\crv{(44,90)};
(48,98)*{^N};
(52,50)*{{\bf =}};
(67,0);
(109,0)**@{.};
(67,20);
(109,20)**@{.};
(67,40);
(109,40)**@{.};
(67,60);
(109,60)**@{.};
(67,80);
(109,80)**@{.};
(67,100);
(109,100)**@{.};
(67,0);
(67,100)**@{--};
(109,0);
(109,100)**@{--};
(60,10)*{\tau^k_k},
(74,2)*{^1},
(109,17);
(74,0)*{\bullet}**@{-};
(67,17);
(74,20)*{\bullet}**@{-},
(102,2)*{^k},
(109,3);
(102,0)*{\bullet}**@{-};
(67,3);
(102,20)*{\bullet}**@{-},
(60,30)*{\beta^k_t},
(70,40)*{\bullet};
(76,40)*{\bullet}**\crv{(73,34)};
(76,21)*{^1},
(74,20);
(82,40)*{\bullet}**\crv{(75,35),(80,37)};
(93,40)*{\bullet};
(100,40)*{\bullet}**\crv{(97,34)};
(104,21)*{^k};
(102,20);
(106,40)*{\bullet}**\crv{(102,30)};
(60,50)*{\nu^*},
(72,41)*{^1},
(70,40);
(72,60)*{\bullet}**\crv{(70,50)};
(78,41)*{^2},
(76,40);
(80,60)*{\bullet}**\crv{(78,55)};
(84,41)*{^3},
(82,40);
(93,40)**\crv{(87,47)};
(103,41)*{^{N-1}};
(100,40);
(96,60)*{\bullet}**\crv{(99,52)};
(107,41)*{^N};
(106,40);
(104,60)*{\bullet}**\crv{(105,52)};
(60,70)*{{\bf 1}_r},
(72,60);
(72,80)*{\bullet}**@{-};
(73,78)*{^1};
(80,60);
(80,80)*{\bullet}**@{-};
(81,78)*{^2};
(96,60);
(96,80)*{\bullet}**@{-};
(97,78)*{^{r-1}};
(104,60);
(104,80)*{\bullet}**@{-};
(105,78)*{^r};
(73,58)*{^1};
(81,58)*{^2};
(97,58)*{^{r-1}};
(106,58)*{^r};
(60,90)*{\mu},
(72,80);
(70,100)*{\bullet}**\crv{(72,90)};
(71,98)*{^1};
(75,98)*{^2};
(76,100)*{\bullet};
(82,100)*{\bullet}**\crv{(79,94)};
(83,98)*{^3};
(80,80);
(93,100)*{\bullet}**\crv{(82,90)};
(96,80);
(100,100)*{\bullet}**\crv{(97,90)};
(100,98)*{^{N-1}};
(104,80);
(106,100)*{\bullet}**\crv{(104,90)};
(108,98)*{^N};
\end{xy}
\]

When we apply either $\tau^r_r$ between $\mu$ and $\nu^*$ or
$\tau^k_k$ below $\beta^k_t$ we get the same result: every through
line of the resulting diagram makes one more turn around the
cylinder.

Then for $k \in \poset r$ we have
$$
\pi_k(a_r) = \sum_{i,j=1}^{d_r} \pi_k \left( \beta^r_i P_{r-2}(\tau^r_r)
C_{i,j}(\tau_r)(\beta^r_j)^* \right) =  \sum_{i,j=1}^{d_r}
\pi_k \left (\beta^r_i C_{i,j}(\tau_r)(\beta^r_j)^* \right)P_{r-2}(x^k)
$$
$$
= P_{r-2}(x^k) \pi_k \left( \sum_{i,j=1}^{d_r}
\beta^r_i C_{i,j}(\tau_r)(\beta^r_j)^*\right).
$$
So if $0 < k < r$ then $\pi_k(a_k)$ is divisible by $P_k(x^k)$ because
$P_{r-2}(x^k)$ is divisible by $P_k(x^k)$. For $k=0$ one has $\pi_0(a_r) =
\zero$ because $P_{r-2}(1)=0$. The theorem is proved.
\end{proof}

Note that Proposition~\ref{centralizer} can be deduced
independently from Theorem~\ref{ideal_inside}.

\begin{proof} Consider the matrix $A = \prod_{k \in \posetN}A_k$ commuting with
$\TL$. Let $k>0$. Then $A_k$ commutes with $\matr {d_k}{\Cox}P_k(x^k)$. Consider
$z \in \Cot \backslash \text{ roots of $P_k(x^k)$}$. For every matrix $B
\in \matr{d_k}{\Co}$ there exist a matrix $C \in \matr {d_k}{\Cox}P_k(x^k)$
such that $C(z) = B$. Thus $A_k(z)$ commutes with $\matr{d_k}{\Co}$, hence it
is diagonal. All nondiagonal entries of $A_k$ equal zero in infinitely many
points of $\Co$, so they are zero everywhere on $\Co$. All diagonal entries of
$A_k$ coincide in infinitely many points, so they coincide as elements
of $\Cox$. We proceed analogously for $k=0$.
\end{proof}

\section{Geometry of the center}

In this section we consider the center $\centr{\TL}$ and prove
Corollary~\ref{variety_decomposition}. We identify $\TL$ with
the subalgebra of $\prod_{k \in \posetN}\matr{d_k}{\Cox}$.

Put
$$
\idealN = \begin{cases}
\prod_{k \in \posetN} P_k(x^k)\Cox &\hbox{if $N$ is odd},\\
P_0(x) \CoxS  \times \prod_{k \in \posetN,k>0} P_k(x^k)\Cox &\hbox{if $N$ is
even},\\
\end{cases}
$$
$$
\overN = \begin{cases}
\prod_{k \in \posetN} \Cox &\hbox{if $N$ is odd},\\
\CoxS  \times \prod_{k \in \posetN,k>0} \Cox &\hbox{if $N$ is
even}.\\
\end{cases}
$$

Then
\begin{equation}\label{center_inclusion}
\idealN \subset \centr{\TL} \subset
\overN
\end{equation}
by Theorem~\ref{ideal_inside} and Proposition~\ref{centralizer}.

\begin{proposition}\label{generators}
The center $\centr{\TL}$ is a finitely generated algebra without nilpotent
elements.
\end{proposition}

To prove the proposition we need the following facts.

\begin{lemma}\label{f_g_polynomial}
Let $P(x) \in \Co[x]$. The ideal $I_P = P(x)\Co[x]$ of $\Co[x]$ is finitely
generated as algebra.
\end{lemma}
\begin{proof}
Let $n = \deg P$. If $n=0$ then $I_P=\Co[x]$ and elements $1, x$ are the
generators. If $n>0$ one can take $Y_0 = P(x), \dots, Y_{n-1} = x^{n-1}P(x)$ as
a set of generators.
\end{proof}

\begin{lemma}\label{f_g_L_polynomial}
Let $P(x) \in \Cox$. The ideal $I_P = P(x)\Cox$ of $\Cox$ is finitely generated
as algebra.
\end{lemma}
\begin{proof}
If $P(x)$ is invertible in $\Cox$ then $I_P = \Cox$ and generators are $x$ and
$\frac 1x$. Recall that the only invertible elements in $\Cox$ are $\alpha x^m$
for some $\alpha \in \Cot, m \in \Z$. If $P(x)$ is equal to $x+\alpha$ for some
$\alpha \in \Cot$ up to multiplication by invertible element of $\Cox$, then
the generators are $x+\alpha$ and $\frac 1{\alpha} + \frac 1x$. Otherwise we
can suppose that $P(x) = \alpha_n x^n+\dots+\alpha_{-m}\frac 1{x^m}$ for some
$n,m \ge 1$ and $\alpha_n,\alpha_{-m} \in \Cot$. Indeed, one can make $P(x)$ to
be of this form multiplying $P$ by some invertible element of $\Cox$. Then
$Y_{-m} = x^{-m}P(x), \dots, Y_{n-2}=x^{n-2}P(x), Y_{n-1} = x^{n-1}P(x) \in I_P$
is the set of generators.
\end{proof}

\begin{proof}[Proof of Proposition~\ref{generators}]
There are no nilpotent elements in $\centr{\TL}$ as it is included into algebra
without nilpotent elements by~(\ref{center_inclusion}).

$\idealN$ is finitely generated algebra as it is a direct product of finitely
generated algebras (Lemma~\ref{f_g_polynomial} and
Lemma~\ref{f_g_L_polynomial}).

$\fact{\centr{\TL}}{\idealN}$ is finite dimensional. Indeed,
$\fact{\centr{\TL}}{\idealN} \subset \fact{\overN}{\idealN}$ and
$$
\fact{\overN}{\idealN} \cong
\begin{cases}
\prod_{k \in \posetN}{\left(\fact{\Cox}{I_{P_k}}\right)} &\hbox{if $N$ is
odd},\\
\left(\fact{\Co[x]}{I_{\widetilde{P_0}}}\right)\times \prod_{k \in \posetN,k>0}
{\left(\fact{\Cox}{I_{P_k}}\right)} &\hbox{if $N$ is
even},\\
\end{cases}
$$
where $\widetilde{P_0} \in \Co[x]$ is such that
$P_0(x)=\widetilde{P_0}(x+\frac 1 x)$.
The algebra on the right is of finite dimension.

Let $s^1,\dots,s^L \in \centr{\TL}$ are such that their
images under factorization modulo $\idealN$ form a linear basis. Then
$\centr{\TL}$ is generated by generators of $\idealN$ together with
$s^1,\dots,s^L$.
\end{proof}

Any finitely generated commutative algebra with $\un$ over an algebraically
closed field and without nilpotent elements is a coordinate ring of a
certain affine algebraic variety (it is an exercise in \cite{Hartshorne}). Let
us denote this variety for $\centr {\TL}$ by $\varietyN$.

We are going to describe $\varietyN$ in this section. We  will find what are its
irreducible components and construct a desingularization of $\varietyN$.

Consider the sets $S_k = \Cot$ for $k>0$ and $S_0 = \Co$. Consider the map
\begin{equation}\label{psi_map}
\Psi: \bigsqcup_{k \in \posetN} S_k \To \varietyN
\end{equation}
defined by the following. Any element $f \in \centr{\TL}$ is a set of Laurent
polynomials $\{f_k; k \in \posetN \}$. For $z \in S_k$ consider a
homomorphism $\psi_{k,z}: \centr{\TL} \To \Co$ defined by
$$
\psi_{k,z}(f) = \begin{cases}
f_k(z) &\hbox{if $k>0$},\\
f_0(t) \hbox{ where $t + \frac 1t = z$} &\hbox{if $k=0$}.\\
\end{cases}
$$
Then $\kernel{\psi_{k,z}}$ is a maximal ideal in $\centr{\TL}$,
so it defines a unique point $\Psi(z) \in \varietyN$.

\begin{theorem}\label{desingularization}
$\image{\Psi} = \varietyN$ and $\Psi$ is one-to-one on the set
$$
\begin{cases}
\bigsqcup_{k \in \posetN} S_k \backslash \left\{ \hbox{roots of } P_k(x^k)
\right\} &\hbox{if $N$ is odd}\\
S_0 \backslash \left\{ t+\frac 1t | P_0(t)=0 \right\} \bigsqcup \bigsqcup_{k
\in \posetN, k>0} S_k \backslash \left\{ \hbox{roots of } P_k(x^k)
\right\} &\hbox{if $N$ is even}\\
\end{cases}
$$
\end{theorem}

To prove this statement we need the following fact.

\begin{lemma}\label{curve}
Let $P_1,\dots,P_L \in \Co[x]$. Consider the map $\phi: \Co \To \Co^L$ given
by $\phi(x) = \left( P_1(x),\dots,P_L(x) \right)$. Then $\image{\phi}$ is
algebraically closed.
\end{lemma}
\begin{proof}
Let $m = max_{1 \le i \le L} \deg{P_i}$. Consider the map of projective spaces
$\hat{\phi}: \Proj^1 \To \Proj^L$ given in homogenous coordinates as
$$
\hat{\phi}([X_0,X_1]) = [X_0^m,\hat{P}_1(X_0,X_1),\dots,\hat{P}_L(X_0,X_1)]
$$
where for the polynomial $Q(x)=\sum_k a_k x^k$ we denote
$\hat{Q}(X_0,X_1) = \sum_k a_k X_1^k X_0^{m-k}$. This map is
correctly defined: for $X_0 \ne 0$ one has nonzero first
coordinate $X_0^m \ne 0$; otherwise $X_1 \ne 0$ and there is a
number $k$ with $\deg P_k = m$, so $\hat{P_k}(0,X_1) \ne 0$. So
$\hat{\phi}$ is regular map of projective variety to projective
variety, thus $\image{\hat{\phi}}$ is closed. There is a
canonically defined open set $U_0 = \{ Y_0 \ne 0 \} \subset
\Proj^L$ isomorphic to $\Co^L$, and the intersection of any closed
projective variety with $U_0$ is isomorphic to a closed affine
variety. Isomorphism is given by
$$
[Y_0,Y_1,\dots,Y_L] \To (\frac {Y_1}{Y_0},\dots,\frac{Y_L}{Y_0}).
$$
Then $\image{\hat{\phi}} \cap U_0$ is isomorphic to $\image{\phi}$. Lemma is
proved.
\end{proof}

Similar statement is true for Laurent polynomials but under additional
assumptions. Let $P(x) = \sum_{i=m}^{n}a_ix^i$ where $a_m,a_n \ne 0$ be Laurent
polynomial. We denote $mindeg P = m$ and $maxdeg P = n$.

\begin{lemma}\label{L_curve}
Let $P_1,\dots,P_L \in \Cox$. Consider the map $\phi: \Cot \To \Co^L$ given
by $\phi(x) = \left( P_1(x),\dots,P_L(x) \right)$. Suppose there exist $m > 0$
such that
$$
\underset{1 \le i \le L} {min} mindeg{P_i} = -m, \; \underset{1 \le i \le
L} {max} maxdeg{P_i} = m
$$
and there exist $i \ne j$ such that $mindeg P_i = -m$, $maxdeg P_i < m$ and
$mindeg P_j > -m$, $maxdeg P_j = m$. Then $\image{\phi}$ is algebraically
closed.
\end{lemma}
\begin{proof}
Consider the map $\hat{\phi}: \Proj^2 \To \Proj^L$ given in
homogenous coordinates as
$$
\hat{\phi}([X_0,X_1,X_2]) =
[X_0^m,\hat{P_1}(X_0,X_1,X_2),\dots,\hat{P_L}(X_0,X_1,X_2)]
$$
where for $Q(x) = c + \sum_{i>0}a_ix^i + \sum_{j>0}b_i \frac 1{x^j}$ we define
$$
\hat{Q}[X_0,X_1,X_2] = c X_0^m + \sum_{i>0}a_iX_1^iX_0^{m-i}
+ \sum_{j>0}b_i X_2^jX_0^{m-j}.
$$
It is again correctly defined because either
$\hat{P_i}(0,X_1,X_2) \ne 0$ or $\hat{P_j}(0,X_1,X_2) \ne 0$. Then
$\image{\phi}$ is isomorphic to $U_0 \cap \image{
\hat{\phi}\left(\{X_0^2 = X_1 X_2\}\right)}$, so it is closed.
\end{proof}

Consider an algebra $\subalgN = \Co \un + \idealN$. $\subalgN$ is obviously
finitely generated (see arguments in proof of Proposition~\ref{generators})
and has no nilpotent elements. So $\subalgN$ is a coordinate ring of some affine
algebraic variety, let us denote this variety by $\varietyI$. Then we can
construct the map
$$
\Phi: \bigsqcup_{k \in \posetN} S_k \To \varietyI
$$
in the same way as $\Psi$ for $\varietyN$.

\begin{proof}[Proof of Theorem~\ref{desingularization}]
First, we choose useful coordinates on variety $\varietyN$. Let $s^1,\dots,s^L
\in \centr{\TL}$ be such that form a linear basis in
$\fact{\centr{\TL}}{\idealN}$ after factorization modulo $\idealN$. Recall
that for $f \in \centr{\TL}$ we denote by $\{f_k; k \in \posetN \}$ its components which are Laurent polynomials.

Due to Lemma~\ref{f_g_polynomial} and Lemma~\ref{f_g_L_polynomial} one can
choose generators $Y^k_1 = x^{m_k}P_k(x^k),Y^k_2 = x^{m_k+1}P_k(x^k), \dots,
Y^k_{L_k} = x^{n_k}P_k(x^k)$ in $P_k(x^k)\Cox$ for $k>0$ and
$Y^0_1 = P_0(x),Y^0_2 = (x+\frac1x)P_0(x), \dots,
Y^0_{L_0} = (x+\frac1x)^{L_0-1}P_0(x)$ in $P_0(x)\CoxS$. Then for each $k>0$ we
can make $m_k$ smaller and $n_k$ large for the $k^{th}$ set of generators so that
Laurent polynomials
$$
Y^k_1,\dots,Y^k_{L_k},s^1_k,\dots,s^L_k
$$
satisfy conditions of Lemma~\ref{L_curve}. So, we suppose this condition to be
satisfied.

If $N$ is even, let us consider such polynomials
$\widetilde{Y^0_1},\dots,\widetilde{Y^0_{L_0}},\widetilde{s^1_0},\dots,
\widetilde{s^L_0} \in \Co[x]$ that $Y^0_i(x)=\widetilde{Y^0_i}(x+\frac1x)$
and $s^i_0(x)=\widetilde{s^i_0}(x+\frac1x)$.

So $\left\{ Y^k_i,s^j | k\in\posetN, 1\le i \le L_k, 1\le j \le L \right\}$ can
be considered as a coordinates on $\varietyN$ because they are generators of
$\centr{\TL}$. Then
$$
\Psi(z) = \begin{cases}
\left(0,\dots,0,Y^k_1(z),\dots,Y^k_{L_k}(z),0,\dots,0,s^1_k(z),\dots,s^L_k(z)
\right) &\hbox{ if $z \in S_k$, $k>0$},\\
\left(\widetilde{Y^0_1}(z),\dots,\widetilde{Y^0_{L_0}}(z),0,\dots,0,0,\dots,0,
\widetilde{s^1_0}(z ) ,
\dots,\widetilde{s^L_0}( z) \right) &\hbox{
if $z \in S_0$.}\\
\end{cases}
$$
Then $\image{\Psi}$ is dense in $\varietyN$ because there is no nonzero element
$f \in \centr{\TL}$ which is zero on $\image{\Psi}$, i.e. $f_k(z) = 0$
for each $k \in \posetN$, $z \in \Cot$. But $\image{\Psi} = \bigcup_{k \in
\posetN}\Psi(S_k)$, and each $\Psi(S_k)$ is closed by Lemma~\ref{L_curve} for
$k>0$ and by Lemma~\ref{curve} for $k=0$. So, $\varietyN = \image{\Psi}$.

We consider $\left\{ Y^k_i | k\in\posetN, 1\le i \le L_k \right\}$ as a
coordinates on $\varietyI$ because they are generators of $\idealN$. Then
$\Phi$ can be written in this coordinates in the similar way as $\Psi$, and
one has $\image{\Phi} = \varietyI$ again. Now $\varietyN \subset
\Co^{\sum_k{L_k}+L}$, $\varietyI \subset \Co^{\sum_k{L_k}}$, and
$$
\Phi = \pi \circ \Psi,
$$
where $\pi$ is a projection in $\Co^{\sum_k{L_k}+L}$ on first $\sum_k{L_k}$
coordinates. Consider $z \in S_k$. If $P_k(z^k)=0$ ($P_0(t) = 0$ with $t +
\frac1t = z$) then $\Phi(z) = (0,\dots,0)$. Otherwise $\Phi(z) \ne (0,\dots,0)$
and $z$ and $k$ both can be calculated from $\Phi(z)$ in obvious way. So,
$\Phi$ is one-to-one on the set $\bigsqcup_{k\in\posetN}S_k \Big\backslash
\Phi^{-1}\left((0,\dots,0)\right)$. Then $\Psi$ is also one-to-one on this set.
\end{proof}

So, $\varietyI$ is a bouquet of $\left[\frac n2\right]+1$ irreducible
curves, it has exactly one singular point.

\begin{proof}[Proof of Corollary~\ref{variety_decomposition}]
Denote $\nu_k = \Psi(S_k)$. This is an irreducible curve with finite number of singular points (points of self-intersection) by Theorem~\ref{desingularization}. Then $\varietyN = \cup_{k \in \posetN}{\nu_k}$ and $\nu_k \cap \nu_l$ is a finite
set due to Theorem~\ref{desingularization} again.
\end{proof}

The desingularization $\Psi$ constructed above gives us a way to
give a proper description of $\varietyN$. For those one should
calculate points on sheets $S_k$ correspondent to the points of
intersections of the curves $\nu_k = \Psi(S_k)$. We deduce from
Theorem~\ref{desingularization} that all intersection
(self-intersection) points are among $\Phi^{-1}(0,\dots,0)$. If
one finds some successive family of central elements in $\TL$ it
will separate some of those points one from each other. Following
proposition gives us one more necessary condition for the
intersections.

\begin{proposition}\label{necessary}
Let $x \in S_k$, $y \in S_m$. Then $\Psi(x) = \Psi(y)$
only if

(i) $P_k(x^k)=0$ (or $x = z + \frac1z$ with $P_0(z)=0$ if $k = 0$)

(ii) $P_m(y^m)=0$ (or $y = z + \frac1z$ with $P_0(z)=0$ if $m = 0$)

(iii) $x^k = y^m$
\end{proposition}
\begin{proof}
Two first statements are due to Theorem~\ref{desingularization}. Let us prove
the third one.

Obviously $\{F(\tau_n^n); F \in \Cox \} \subset \centr{\TL}$. Then $\Psi(x) =
\Psi(y)$ only if $F(\tau_n^n)_k(x) = F(\tau_n^n)_m(y)$ for each $F$. But
$F(\tau_n^n)_k(x) = F(x^k)$ and $F(\tau_n^n)_m(y) = F(y^m)$. So, $x^k =
y^m$.
\end{proof}

\section{Example for $\TLN 3$}

Let us consider $\TLN 3$ now. $\poset 3 = \{1, 3\}$, $d_3=1$ ($d_N=1$ because
the only standard diagram in $\diag NN$ is $\un_N$). Following diagrams are
the standard ones in $\diag13$
\[
\begin{xy}
<0cm,0.5cm>;
<1.6cm,0.5cm>**@{.};
<1.6cm,1.5cm>**@{--};
<0cm,1.5cm>**@{.};
<0cm,0.5cm>**@{--},
(3,15)*{\bullet};
(8,15)*{\bullet}**\crv{(5,11)},
(13,15)*{\bullet};
(8,5)*{\bullet}**\crv{(8,10)};
(3,18)*{1},
(8,18)*{2},
(13,18)*{3},
(7,8)*{1},
(8,0)*{\beta_1}
\end{xy}\;\;\;
\begin{xy}
<0cm,0.5cm>;
<1.6cm,0.5cm>**@{.};
<1.6cm,1.5cm>**@{--};
<0cm,1.5cm>**@{.};
<0cm,0.5cm>**@{--},
(8,15)*{\bullet};
(13,15)*{\bullet}**\crv{(10,11)},
(3,15)*{\bullet};
(8,5)*{\bullet}**\crv{(8,10)};
(3,18)*{1},
(8,18)*{2},
(13,18)*{3},
(9,8)*{1},
(8,0)*{\beta_2}
\end{xy}\;\;\;\begin{xy}
<0cm,0.5cm>;
<1.6cm,0.5cm>**@{.};
<1.6cm,1.5cm>**@{--};
<0cm,1.5cm>**@{.};
<0cm,0.5cm>**@{--},
(0,13);
(3,15)*{\bullet}**\crv{(1,13)};
(16,13);
(13,15)*{\bullet}**\crv{(15,13)},
(8,15)*{\bullet};
(8,5)*{\bullet}**\crv{(8,10)};
(3,18)*{1},
(8,18)*{2},
(13,18)*{3},
(6,8)*{1},
(8,0)*{\beta_3}
\end{xy}
\]
where we omit upper indices ($\beta_i$ instead of $\beta^1_i$). So, $d_1=3$
and
$$
\TLN 3 = \left\{ S(\tau_3) + \sum_{i,j=1}^{3}\beta_i F_{i,j}(\tau_1)\beta_j^*
\Big | S, F_{i,j} \in \Cox \right\} \subset \matr 1{\Cox} \times \matr3{\Cox}.
$$


\begin{proposition}\label{centr_TL3}
$\centr{\TLN 3} \cong \left\{ S \times T \in \Cox^2 \Big| S(x) - T(x^3)\;
\vdots \; x^2 + \frac1{x^2}+\delta \right\}$
\end{proposition}
\begin{proof}

For $P \in \Cox$ denote
$$
\hat P = \begin{pmatrix}
    P_0(x)&\frac1x P_2(x)& P_1(x)\\
    x P_1(x)&P_0(x)&P_2(x)\\
    P_2(x)&P_1(x)&P_0(x)\\
    \end{pmatrix},
$$
where $P(x) = P_0(x^3)+xP_1(x^3)+\frac1xP_2(x^3)$. Then
$$
\hat P + \hat Q = \widehat{P + Q}, \;\hat P \hat Q = \widehat{P Q}
$$
and if $\hat P A = A \hat P$ is a scalar matrix then $A = \hat Q$ for some $Q$.

Let $a = S(\tau_3) + \sum_{i,j=1}^{3}\beta_i F_{i,j}(\tau_1)\beta_j^* \in \TLN
3$. Then $a$ acts on $M_3 \times M_1$ by matrix
$$
S(x) \times \begin{pmatrix}
    \hat S + F R_1
    \end{pmatrix}
$$
where $F = \begin{pmatrix}F_{i,j}\end{pmatrix}$ and
$$
R_1 = \left( (\beta_i)^*\beta_j \right)_{i,j=1}^3 = \begin{pmatrix}
\delta & 1 & \frac1x \\
1 & \delta & x \\
x & \frac1x & \delta\\
\end{pmatrix}
$$
was defined by~(\ref{mod_k_matrix}). Indeed, consider
the action of $\tau_3$ on $M_1$:
$$
\tau_3 \beta_1 = \beta_2 \tau_1 = \beta_2 x,\;
\tau_3 \beta_2 = \beta_3,\;
\tau_3 \beta_3 = \beta_1.
$$
Then
$$
\tau_3= x \times \begin{pmatrix}
    0&0&1\\
    x&0&0\\
    0&1&0\\
\end{pmatrix}
,\;\;
\tau_3^2= x^2 \times \begin{pmatrix}
    0&1&0\\
    0&0&x\\
    x&0&0\\
\end{pmatrix}
,\;\;
\tau_3^3= x^3 \times \begin{pmatrix}
    x&0&0\\
    0&x&0\\
    0&0&x\\
    \end{pmatrix},
$$
thus $S(\tau_3) = S(x) \times \hat S$. Also $\sum_{i,j=1}^{3}\beta_i
F_{i,j}(\tau_1)\beta_j^* = \zero \times F R_1$ due to~(\ref{action_matrix}).

$a \in \centr{\TLN 3}$ if and only if $a = S(x) \times T(x)\un_3$
for some $T \in \Cox$ due to Proposition~\ref{centralizer}. Hence,
$\hat S + F R_1 = T \un_3$. Then $F = \hat V$ for some $V \in \Cox$. Indeed,
$F \det R_1 = \widehat{\left(T(x^3)-S(x)\right)} N_1$ where $N_1$ was defined
by~(\ref{inverse}) such that $R_1 N_1 = (\det R_1) \un_3$. Note that $R_1 =
\widehat{x^2 + \frac1{x^2}+\delta}$. Then $N_1 = \hat Q$ for some $Q \in \Cox$,
and $F \det R_1 = \widehat{\left(T(x^3)-S(x)\right)Q(x)}$.

Now, $S(x) + V(x)(x^2 + \frac1{x^2}+\delta) = T(x^3)$. Obviously, one can take
an arbitrary $V \in \Cox$ and put $F = \hat V$.
\end{proof}

We have an explicit description of $\centr{\TLN 3}$. Recall the
desingularization
$$
\Psi: \Cot \sqcup \Cot \To \variety 3
$$
of corresponding affine algebraic variety introduced by~(\ref{psi_map}). We take
$z$ as a coordinate on $S_3 = \Cot$ and $w$ as a coordinate on $S_1 =
\Cot$.

\begin{proposition}
Consider the points $z_1 = \sqrt{q}, z_2 = -\sqrt{q}, z_3 = \frac
1{\sqrt{q}}, z_4 = -\frac 1{\sqrt{q}}$, put $w_i = z_i^3$ for $i=1,2,3,4$.

Then $\Psi(z_i)=\Psi(w_i)$ for $i=1,2,3,4$ and $\Psi$ maps
$\Cot \setminus \{z_i|i=1,2,3,4\} \sqcup \Cot \setminus \{w_i|i=1,2,3,4\}$ onto
$\variety 3 \setminus  \{ \Psi(z_i) | i=1,2,3,4 \}$ bijectively.
\end{proposition}

\begin{proof}
$\Psi(w_1)=\Psi(w_2)$ only if $w_1=w_2$ due to
Proposition~\ref{necessary}~(iii).

Put $P(x)= x^2 + \frac1{x^2}-q^3-\frac1{q^3}$ where $\delta = -q-\frac1q$. Now
$P_3(x)=P_1(x)=P(x)$, and for $w$ such that $P(w)=0$ the map $\Psi$ may glue
some cubic roots of $w$ either with $w$ or with each other. We will consider
these two possibilities.

(I) Let $\Psi(z)=\Psi(w)$. This is equivalent to
$$
\left\{ S \times T \in \centr{\TLN 3} \Big| S(z) = 0 \right\} =
\left\{ S \times T \in \centr{\TLN 3} \Big| T(w) = 0 \right\}
$$
by definition of the map $\Psi$. Recalling Proposition~\ref{centr_TL3} the last
is equivalent to
$$
\left\{ (V,T) \in \Cox^2 \Big| T(z^3) + V(z)( z^2+\frac1{z^2}+\delta ) = 0
\right\}
$$
$$
= \left\{ (V,T) \in \Cox^2 \Big| T(w)=0 \right\},
$$
or
\begin{equation}\label{glueing_condition}
z^2+\frac1{z^2}+\delta = 0
\end{equation}
since $z^3 = w$.

Conversely, if~(\ref{glueing_condition}) holds then $\Psi(z)=\Psi(w)$ where
$w=z^3$.

(II) Let $\Psi(z)=\Psi(x)$ for some unequal $z,x \in S_3$. Then
$z^3=x^3$ by Proposition~\ref{necessary}~(iii), denote $w = z^3$.
If $z$ and $x$ are roots of~(\ref{glueing_condition}) both then
$\Psi(z)=\Psi(x)=\Psi(w)$ by part (I). Suppose $z$ is not a root
of~(\ref{glueing_condition}), put $a = z^2+\frac1{z^2}+\delta \ne
0$. Then $\Psi(z)=\Psi(x)$ implies
$$
\left\{ (V,T) \in \Cox^2 \Big| V(z) = -\frac{T(w)}{a} \right\}
$$
$$
= \left\{ (V,T) \in \Cox^2 \Big| V(x)( x^2+\frac1{x^2}+\delta ) = -T(w)
\right\}
$$
what is obviously a contradiction.
\end{proof}

So, for generic $q \in \Cot$ variety $\variety 3$ looks like two copies of
$\Cot$ glued in four points ($z_i = w_i$). If $q^6=1$ some of those points
coincide, so the number of singular points will be less than~4. Especially:

$q=1$. Then we have two points $z=-1,1$ glued to $w=-1,1$ correspondingly.

$q=\frac{1 \pm 3 \sqrt{-1}}2$. Then we have two points
$z=\alpha,-\frac1{\alpha}$ glued to $w=i$ and two points
$z=-\alpha,\frac1{\alpha}$ glued to $w=-i$
correspondingly where $\alpha = \cos{\frac{\pi}6} +
\sqrt{-1}\sin{\frac{\pi}6}$.

$q=\frac{-1 \pm 3 \sqrt{-1}}2$. Then we have two points
$z=\alpha,\frac1{\alpha}$ glued to $w=-1$ and two points
$z=-\alpha,-\frac1{\alpha}$ glued to $w=1$ correspondingly
where $\alpha = \frac{1 + 3 \sqrt{-1}}2$.

$q=-1$. Then we have two points $z=-\sqrt{-1},\sqrt{-1}$ glued to
$w=\sqrt{-1},-\sqrt{-1}$ correspondingly.

Case $q=-1$ ($\delta = 2$) is when the corresponding Hecke algebra is a group
algebra of the affine Coxeter group.

%

\end{document}